# BASIC VACCINATION CONTROL TECHNIQUES FOR A GENERAL TRUE -MASS ACTION SEIR MODEL

M. De la Sen and S. Alonso-Quesada Department of Electricity and Electronics. **Faculty of Science and Technology**. University of the Basque Country. Campus of Leioa, Aptdo. 644 - Bilbao. SPAIN.

**Abstract**: This paper presents several simple linear vaccination-based control strategies for a SEIR (susceptible plus infected plus infectious plus removed populations) propagation disease model. The model takes into account the total population amounts as a refrain for the illness transmission since its increase makes more difficult contacts among susceptible and infected. The vaccination control objective is the asymptotically tracking of the removed-by-immunity population to the total population while achieving simultaneously that the remaining populations (i.e. susceptible plus infected plus infectious) tend asymptotically to zero.

## 1. INTRODUCTION

Important control problems nowadays related to Life Sciences are the control of ecological models like, for instance, those of population evolution (Beverton-Holt model, Hassell model, Ricker model etc.) via the online adjustment of the species environment carrying capacity, that of the population growth or that of the regulated harvesting quota [1-5], as well as those of disease propagation via vaccination control. Several variants and generalizations of the Beverton-Holt model (standard time-invariant, time-varying parameterized, generalized model or modified generalized model) have been investigated at the levels of stability, cycle-oscillatory behaviour, permanence and control through the manipulation of the carrying capacity, [1, 4]. The design of related control actions has been proved to be relevant at the levels, for instance, of aquaculture exploitation or plague fighting. At the same time, properties of the discrete Beverton-Holt equation solutions when the species environment carrying capacity and/or the intrinsic growth rate of the species are periodic functions of time have been studied, [5]. Such a particular case is of interest since periodic fluctuations are quite common in Biology and Ecology. For instance, the habitat resource availability, temperature, humidity and so on, which influence the environment carrying capacity, experiment fluctuation during the year. In the same way, the birth and/or the survivorship rates of the species may suffer fluctuations along the year, so a periodical intrinsic growth rate is an appropriate choice in such ecology models. Furthermore, this kind of difference equations for modelling real life situations in population biology, ecology, economics and so on from a control theory point of view has been recently dealt with in a generalized way in [6]. The boundedness properties of the positive solutions of this type of nonlinear difference equations have been established in such a paper.

On the other hand, the literature about epidemic mathematical models is exhaustive in many books and papers (see, for instance, [7-16]). The sets of models include the most basic ones:

- SI-models, where only susceptible and infected populations are assumed,

- SIR models, which include susceptible plus infected plus removed-by-immunity populations, and
- SEIR-models, where the infected populations is split into two ones (namely, the "infected" which incubate the disease but do not still have any disease symptoms and the "infectious" or "infective" which do have the external disease symptoms).

Those models have also two major variants, namely, the so-called "pseudo-mass action models", where the total population is not taken into account as a relevant disease contagious factor and the so-called "true-mass action models", where the total population is more realistically considered as an inverse factor of the disease transmission rates. There are many variants of the above models, for instance, including vaccination of different kinds: constant [9], impulsive [13], discrete-time etc., incorporating point or distributed delays [13, 14], oscillatory behaviors [15] and so on. On the other hand, variants of such models become considerably simpler for the illness transmission among plants [7, 8].

In this paper, two continuous-time vaccination control strategies are given for a SEIR epidemic model. One of them takes directly the susceptible population to design the vaccination strategy while the second one uses the total and removed-by-immunity populations for such a purpose. It is assumed that the total population remains constant through time, so that the illness transmission is not critical, and the SEIR-model is of the above mentioned true-mass action type. The positivity of the mathematical model, which reflects the real problem at hand, is proved. Such a property is crucial to guarantee the boundedness for all time of all the partial populations.

## 2. SEIR-EPIDEMIC MODEL

Let $S(t)$ be the "susceptible" population of infection, $E(t)$ the "infected" population, $I(t)$ the "infectious" population and $R(t)$ the "removed-by-immunity" (or "immune") population at time $t$. Consider the true-mass action type SEIR epidemic model:

$$\dot{S}(t) = -\mu S(t) + \omega R(t) - \beta \frac{S(t)I(t)}{N} + \mu N (1 - V(t)) \quad (1)$$

$$\dot{E}(t) = \beta \frac{S(t)I(t)}{N} - (\mu + \sigma)E(t) \quad (2)$$

$$\dot{I}(t) = -(\mu + \gamma)I(t) + \sigma E(t) \quad (3)$$

$$\dot{R}(t) = -(\mu + \omega)R(t) + \gamma I(t) + \mu N V(t) \quad (4)$$

subject to initial conditions $S(0) \geq 0$, $E(0) \geq 0$, $I(0) \geq 0$ and $R(0) \geq 0$ under the vaccination constraint $V: \mathbb{R}_{0+} \to \mathbb{R}_{0+}$ where $\mathbb{R}_{0+} \triangleq \{z \in \mathbb{R} \mid z \geq 0\}$. In such a SEIR-model, $N$ is the constant total population, $\mu$ is the rate of deaths from causes unrelated to the infection, $\omega$ is the rate of losing immunity, $\beta$ is the transmission constant (with the total number of infections per unity of time at time $t$ being $\beta \frac{S(t)I(t)}{N}$), $\sigma^{-1}$ and $\gamma^{-1}$ are, respectively, the average durations of the latent and infective periods. All the above parameters are assumed to be nonnegative.

The following elementary result follows from the SEIR mathematical model (1)-(4):

**Assertion 1**. The SEIR model (1)-(4) fulfils the constant population through time constraint, i.e.:

$$N(t) \triangleq S(t) + E(t) + I(t) + R(t) = N(0) = N_0 = N > 0 \qquad (5)$$

irrespective of the vaccination strategy.

*Proof*: It follows immediately by summing-up both sides of (1) to (4) what leads to:

$$\dot{N}(t) = \dot{S}(t) + \dot{E}(t) + \dot{I}(t) + \dot{R}(t) = \mu\left(N(t) - S(t) - E(t) - I(t) - R(t)\right) = 0 \quad \forall t \in \mathbb{R}_{0+} \qquad (6)$$

so that:

$$N(0) = S(0) + E(0) + I(0) + R(0) = N_0 \;\Rightarrow\; N(t) = S(t) + E(t) + I(t) + R(t) = N_0 = N \quad \forall t \in \mathbb{R}_{0+} \quad (7)$$
\*\*\*

**Remark 1**. The ideal vaccination mechanism objective is to reduce to zero the numbers of susceptible, infected and infectious independent of their initial numbers so that the total population becomes equal to the removed-by-immunity population after a certain time. After inspecting (1) and (4), it becomes obvious that the constraint $V: \mathbb{R}_{0+} \to \mathbb{R}_{0+}$ is necessary to decrease the time variation of the susceptible and to increase simultaneously that of the removed by immunity for all time. However, Assertion 1 proves that the constant population through time is independent of the vaccination strategy so that it is independent of the ideal vaccination objective constraint $V: \mathbb{R}_{0+} \to \mathbb{R}_{0+}$ as a result. For instance, in a biological war, the objective would be to increase the numbers of the infected plus the infectious population for all time. For that purpose, the appropriate vaccination strategy would be negative. \*\*\*

The fact that the total population of the SEIR model (1)-(4) remains constant (Assertion 1) makes it both uncontrollable to the origin and unreachable. Such a constraint is atypical in most of control problems since the role of the vaccination is to decrease to zero the numbers of susceptible, infected and infectious to make the removed-by-immunity population to asymptotically converge to the total population. In this context, the following elementary result follows:

**Assertion 2**. The SEIR model (1)-(4) is unreachable and uncontrollable to the origin by using any vaccination strategy.

*Proof*: Proceed by contradiction. Fix any desired final state $x^* \triangleq \begin{bmatrix} S(t^*) & E(t^*) & I(t^*) & R(t^*) \end{bmatrix}^T$ at arbitrary finite time $t = t^*$ fulfilling the constraint $S(t^*) + E(t^*) + I(t^*) + R(t^*) > N$. From Assertion 1, the population remains constant equal to $N$ so that $x^*$ is unreachable at any time $t^*$. Thus, the SEIR model is unreachable. It is always trivially uncontrollable to the origin for arbitrary initial conditions for the total population. ***

*2.1. About the positivity of the SEIR epidemic model*

The vaccination strategy has to be implemented so that the SEIR model be positive in the usual sense that none of the populations, namely, susceptible, infected, infectious and immune be negative at any time. This requirement follows directly from the nature of the problem at hand. The following assumption is made:

**Assumption 1**. The following constraints are assumed on the SEIR-model (1)-(4): $\min\{S(0), I(0), R(0)\} \geq 0$, $E(0) > \frac{\mu+\gamma}{\sigma} I(0)$ and $E(0) < \frac{\beta S(0) I(0)}{(\mu+\sigma)N}$ if $I(0) \neq 0$. ***

**Remark 2**. The physical interpretation of Assumption 1 is that the time origin of interest to fix initial conditions in the SEIR model is the time instant at which the disease starts to be infectious. In this sense, the growing rate of infectious at the time origin must be positive, i.e. $\dot{I}(0) > 0$, even under zero initial condition $I(0) = 0$. In such a particular case, $\dot{E}(0) + \dot{I}(0) < 0$, with $\dot{E}(0) < 0$ from (2) and $\dot{I}(0) > 0$ from (3) since $E(0) > 0$ is derived from Assumption 1. However, if $I(0) \neq 0$, Assumption 1 implies that $\dot{E}(0) + \dot{I}(0) > 0$, with $\dot{E}(0) > 0$ from (2) and $\dot{I}(0) > 0$ from (3). Moreover, from Assumption 1 is derived that:

- $S(0) < N - R(0) - \left(1 + \frac{\mu+\gamma}{\sigma}\right) I(0)$ and

- $S(0) > \frac{(\mu+\gamma)(\mu+\sigma)}{\sigma\beta} N$ if $I(0) \neq 0$

which requires that $\beta > \beta_0 \triangleq (\mu+\gamma)\left(1+\frac{\mu}{\sigma}\right)$ since $S(0) < N$. Such a parametrical condition is of interest even if $I(0) = 0$ in order to make the SEIR model parameters independent of any set of admissible initial conditions. ***

The following result about the positivity of the SEIR model follows:

**Theorem 1**. Assume a vaccination function $V: \mathbb{R}_{0+} \to [0, 1]$ and that the initial conditions satisfy Assumption 1. Then, all the solutions of the SEIR-model (1)-(4) satisfy $S(t), E(t), I(t), R(t) \in [0, N]$ $\forall t \in \mathbb{R}_{0+}$.

***Proof***: The constant population constraint (5) is used in (1), (3) and (4) to eliminate the infected population $E(t)$ leading to:

$$\dot{S}(t) = -(\mu+\alpha)S(t) + \omega R(t) + \left(\alpha - \beta \frac{I(t)}{N}\right)S(t) + \mu N(1-V(t)) \tag{8}$$

$$\dot{I}(t) = -(\mu + \gamma + \sigma)I(t) + \sigma(N - S(t) - R(t)) \quad (9)$$

$$\dot{R}(t) = -(\mu + \omega)R(t) + \gamma I(t) + \mu N V(t) \quad (10)$$

for any given real constant $\alpha \geq (\beta/N)\sup_{t\geq 0}\{I(t)\}$. Such a constraint is guaranteed with $\alpha \geq \alpha_0 = \beta$ if $0 \leq I(t) \leq N$ for all $t \geq 0$. It is possible to rewrite (8)-(10) in a compact form as a dynamic system of state $x(t) = [S(t) \ I(t) \ R(t)]^T$, output $y(t) = S(t) + R(t)$ and whose input is appropriately related to the vaccination function as $u(t) = [1 - V(t) \ V(t)]^T$. This leads to the following set of identities:

$$\dot{x}(t) = A(\alpha)x(t) + \mu N \bar{E}_{13} u(t) + \left(\left[\left(\alpha - \beta \frac{I(t)}{N}\right) E_1 - \sigma E_{13}\right] x(t) + \sigma N e_2\right)$$

$$= A(\alpha)x(t) + \mu N e_3 V(t) + \left(\left(\alpha - \beta \frac{I(t)}{N}\right) E_1 x(t) + \sigma(E(t) + I(t))e_2 + \mu N e_1(1 - V(t))\right) \quad (11)$$

$$y(t) = e_{13}^T x(t)$$

where $e_i \in \mathbb{R}^3$ is the i-th unit Euclidean column vector with its i-th component being equal to one and the other two components being zero, $e_{13} = [1 \ 0 \ 1]^T$ and

$$A(\alpha) \triangleq \begin{bmatrix} -(\mu + \alpha) & 0 & \omega \\ 0 & -(\mu + \gamma + \sigma) & 0 \\ 0 & \gamma & -(\mu + \omega) \end{bmatrix}; E_{13} \triangleq \begin{bmatrix} 0 & 0 & 0 \\ 1 & 0 & 1 \\ 0 & 0 & 0 \end{bmatrix}; \bar{E}_{13} \triangleq [e_1 \ e_3] = \begin{bmatrix} 1 & 0 \\ 0 & 0 \\ 0 & 1 \end{bmatrix}; E_1 \triangleq \begin{bmatrix} 1 & 0 & 0 \\ 0 & 0 & 0 \\ 0 & 0 & 0 \end{bmatrix}$$

(12)

Note that:

**(i)** $A(\alpha)$ is a Metzler matrix [17] for any given $\alpha \in \mathbb{R}_{0+}$.

**(ii)** $\min\{\sigma, \mu\} \geq 0$, $\alpha \geq \beta \frac{I(t)}{N}$, $V(t) \in [0, 1]$ $\forall t \in \mathbb{R}_{0+}$, $e_1, e_2, e_3, e_{13} \in \mathbb{R}^3_{0+}$ and $E_1 \in \mathbb{R}^{3\times 3}_{0+}$.

**(iii)** From Assumption 1 (see also Remark 2), $E(0) + I(0) > 0$ and $\dot{E}(0) + \dot{I}(0) > 0$ if $I(0) > 0$, what also requires that $N \geq S(0) > \frac{(\mu + \sigma)(\mu + \gamma)}{\beta \sigma} N$ and $\beta > (\mu + \gamma)\left(1 + \frac{\mu}{\sigma}\right)$. Then, from continuity of any solution of (1)-(4), it exists $t_1 > 0$ such that $E(t) + I(t) > 0$ $\forall t \in (0, t_1)$. Otherwise, i.e. if $I(0) = 0$, $\dot{E}(0) + \dot{I}(0) < 0$ and $\dot{I}(0) > 0$ implying that $\dot{E}(0) < 0$. Again from continuity arguments, it exists $t_1 > 0$ such that $0 < E(t) < E(0)$ $\forall t \in (0, t_1)$.

Then, one has that, for any admissible initial condition $x(0) = [S(0) \ I(0) \ R(0)]^T$, the unique solution on $[0, t_1)$ of (11) is:

$$x(t) = e^{A(\alpha)t}\left(x(0) + \int_0^t e^{-A(\alpha)\tau} m(\tau) d\tau\right) \in \mathbb{R}^3_+ \quad \forall t \in [0, t_1] \quad (13)$$

since $e^{A(\alpha)t} x(0) \in \mathbb{R}^3_+$ and

$$m(t) = m_1(t) + \sigma\big(E(t) + I(t)\big)e_2 = m_1(t) + \sigma\big(E(t) + e_2^T x(t)\big)e_2 \in \mathbb{R}_+^3 \quad \forall t \in [0, t_1] \tag{14}$$

due to, on one hand:

$$m_1(t) \triangleq \mu N e_3 V(t) + \left(\left(\alpha - \beta\frac{I(t)}{N}\right)E_1 x(t) + \mu N e_1\big(1 - V(t)\big)\right) \in \mathbb{R}_{0+}^3 \quad \forall t \in [0, t_1] \tag{15}$$

and, on the other hand, from (2):

$$E(t) = e^{-(\mu+\sigma)t}\left(E(0) + \frac{\beta}{N}e_1^T\left(\int_0^t e^{(\mu+\sigma)\tau}x(\tau)x^T(\tau)d\tau\right)e_2\right) \in \mathbb{R}_+^3 \quad \forall t \in [0, t_1] \tag{16}$$

Since $x(t_1) \in \mathbb{R}_+^3$ and $e^{-(\mu+\sigma)t}E(0) \in \mathbb{R}_+^3$ $\forall t \in \mathbb{R}_{0+}$, it exists $t_2 > t_1$ such that $E(t) \in \mathbb{R}_+$, $m(t) \in \mathbb{R}_+^3$ and $x(t) \in \mathbb{R}_+^3$ (so that $S(t), I(t), R(t) \in \mathbb{R}_{0+}$) $\forall t \in [0, t_2]$. The above properties extend to $t \in \mathbb{R}_{0+}$ from the structures of (13)-(16). Furthermore, $\lim_{t \to \infty} \inf \{x(t)\} \in \mathbb{R}_{0+}^3$ and $\lim_{t \to \infty} \inf \{E(t)\} \in \mathbb{R}_{0+}$. These relations also imply from (5) that $\max\{S(t), E(t), I(t), R(t)\} \leq N$ $\forall t \in \mathbb{R}_{0+}$. ***

**Remark 3**. Note that the mathematical SEIR-model is not guaranteed to be positive according to Theorem 1 in the sense of [17] since Assumption 1 establishes conditions on the initial conditions. ***

**Corollary 1**. Theorem 1 still holds if $V(t) \in \left[0, \ 1 + \left(\alpha - \beta\frac{I(t)}{N}\right)\frac{S(t)}{\mu N}\right]$ $\forall t \in \mathbb{R}_{0+}$.

*Proof*: It follows from the proof of Theorem 1 since $m_1(t) \in \mathbb{R}_{0+}^3$ $\forall t \in \mathbb{R}_{0+}$ from (15) under this modified vaccination constraint. ***

*2.2. Vaccination-free case equilibrium points and stability*

The equilibrium points $x^* = \begin{bmatrix} S^* & I^* & R^* \end{bmatrix}^T$ of (8)-(10) under identically zero vaccination strategy satisfy the set of constraints:

$$-\mu S^* + \omega R^* - \beta\frac{S^* I^*}{N} + \mu N = 0 \tag{17}$$

$$-(\mu + \gamma + \sigma)I^* + \sigma(N - S^* - R^*) = 0 \tag{18}$$

$$-(\mu + \omega)R^* + \gamma I^* = 0 \tag{19}$$

From (19), it follows that:

$$R^* = \frac{\gamma}{\mu + \omega}I^* \tag{20}$$

By introducing (20) in (18), one obtains that:

$$I^* = \frac{\sigma(\mu + \omega)(N - S^*)}{(\mu + \gamma + \sigma)(\mu + \omega) + \gamma\sigma} \tag{21}$$

By using (20)-(21) in (17), it follows that:

$$\left[\mu + \left(\frac{\omega\gamma}{\mu + \omega} - \frac{\beta}{N}S^*\right)\frac{\sigma(\mu + \omega)}{(\mu + \gamma + \sigma)(\mu + \omega) + \gamma\sigma}\right](N - S^*) = 0 \tag{22}$$

A solution of (22) is $S^* = N$ what leads to the feasible equilibrium point $x_1^* = \begin{bmatrix} N & 0 & 0 \end{bmatrix}^T$ from (20)-(21). Such a point is referred to the *disease-free equilibrium point*. If $S^* \neq N$, another solution of (22) is:

$$S^* = \frac{(\mu+\sigma)(\mu+\gamma)}{\sigma\beta} N \triangleq S_2^* \tag{23}$$

provided that $(\mu+\sigma)(\mu+\gamma) \leq \sigma\beta$ is satisfied by the model parameters. By using (23) in (20)-(21), one obtains that:

$$I^* = \frac{(\mu+\omega)[\sigma\beta-(\mu+\sigma)(\mu+\gamma)]}{\beta[(\mu+\gamma+\sigma)(\mu+\omega)+\gamma\sigma]} N \triangleq I_2^* \quad ; \quad R^* = \frac{\gamma[\sigma\beta-(\mu+\sigma)(\mu+\gamma)]}{\beta[(\mu+\gamma+\sigma)(\mu+\omega)+\gamma\sigma]} N \triangleq R_2^* \tag{24}$$

Then, another feasible equilibrium point is $x_2^* = \begin{bmatrix} S_2^* & I_2^* & R_2^* \end{bmatrix}^T$, which is referred to the *endemic equilibrium point*. Finally, from Assertion 1 the infected population at the equilibrium points $x_1^*$ and $x_2^*$ are, respectively, $E_1^* = 0$ and

$$E_2^* = \frac{(\mu+\gamma)(\mu+\omega)[\sigma\beta-(\mu+\sigma)(\mu+\gamma)]}{\sigma\beta[(\mu+\gamma+\sigma)(\mu+\omega)+\gamma\sigma]} N \tag{25}$$

Note that the condition $(\mu+\sigma)(\mu+\gamma) \leq \sigma\beta$ guarantees both $0 < S_2^* \leq N$, $0 \leq E_2^* < N$, $0 \leq I_2^* < N$ and $0 \leq R_2^* < N$ as the real system requires. Finally, note that in the case that $(\mu+\sigma)(\mu+\gamma) = \sigma\beta$ both feasible equilibrium points degenerate in a unique one since $x_1^* = x_2^*$ is derived.

The following result concerning with the local stability of the vaccination-free SEIR mathematical model around its equilibrium points is proven:

**Theorem 2**.

**(i)** The vaccination-free SEIR model (8)-(10) is locally stable around $x_1^*$ if $\beta < \frac{(\mu+\gamma)(\mu+\sigma)}{\sigma}$ and

**(ii)** it is locally stable around $x_2^*$ if $\left\| \frac{\tilde{p}(s)}{p(s)} \right\|_\infty \triangleq \max_{\omega \in \mathbb{R}_{0+}} \left\{ \left| \frac{\tilde{p}(j\omega)}{p(j\omega)} \right| \right\} < \frac{1}{\beta}$, where $j = \sqrt{-1}$ is the complex unit, $\|\cdot\|_\infty$ is the norm of strictly stable transfer functions in the Hardy space $RH_\infty$ and, $p(s)$ and $\tilde{p}(s)$ are the following polynomials:

$$\begin{aligned} p(s) &= (s+\mu+\omega)\left[(s+\mu)(s+\mu+\sigma+\gamma)+\sigma\gamma\right] \\ \tilde{p}(s) &= \frac{1}{N}\left[((s+\mu+\sigma+\gamma)(s+\mu+\omega)+\sigma\gamma)I_2^* - \sigma(s+\mu+\omega)S_2^*\right] \end{aligned} \tag{26}$$

*Proof*: The linearized model (8)-(10) about its equilibrium points is:

$$\begin{bmatrix} \Delta\dot{S}(t) \\ \Delta\dot{I}(t) \\ \Delta\dot{R}(t) \end{bmatrix} = \begin{bmatrix} -\mu - \frac{\beta}{N}I^* & -\frac{\beta}{N}S^* & \omega \\ -\sigma & -(\mu+\sigma+\gamma) & -\sigma \\ 0 & \gamma & -(\mu+\omega) \end{bmatrix} \begin{bmatrix} \Delta S(t) \\ \Delta I(t) \\ \Delta R(t) \end{bmatrix} \quad (27)$$

Then:

**(i)** At the equilibrium point $x_1^*$, the linearized system (27) becomes:

$$\begin{bmatrix} \Delta\dot{S}(t) \\ \Delta\dot{I}(t) \\ \Delta\dot{R}(t) \end{bmatrix} = \begin{bmatrix} -\mu & -\beta & \omega \\ -\sigma & -(\mu+\sigma+\gamma) & -\sigma \\ 0 & \gamma & -(\mu+\omega) \end{bmatrix} \begin{bmatrix} \Delta S(t) \\ \Delta I(t) \\ \Delta R(t) \end{bmatrix} \quad (28)$$

The characteristic equation of (28) becomes:

$$p_1(s) = (s+\mu+\omega)\big[(s+\mu)(s+\mu+\sigma+\gamma) + \sigma(\gamma-\beta)\big] = 0 \quad (29)$$

The characteristic zeros are $-(\mu+\omega)$ and $-\frac{1}{2}\big(2\mu+\sigma+\gamma \pm \sqrt{(\sigma-\gamma)^2 + 4\sigma\beta}\big)$. As a result, the equilibrium point $x_1^*$ of (27) is locally asymptotically Lyapunov stable if $0 \le \beta < \frac{(\mu+\sigma)(\mu+\gamma)}{\sigma}$ since $\mu+\omega > 0$ by definition of parameters $\mu$ and $\omega$.

**(ii)** For the equilibrium point $x_2^*$, the linearized system (27) has a characteristic equation given by:

$$p_2(s) = Det\left\{\begin{bmatrix} s+\mu+\frac{\beta}{N}I_2^* & \frac{\beta}{N}S_2^* & -\omega \\ \sigma & s+\mu+\sigma+\gamma & \sigma \\ 0 & -\gamma & s+\mu+\omega \end{bmatrix}\right\} = p(s) + \beta\tilde{p}(s) = p(s)\left(1+\beta\frac{\tilde{p}(s)}{p(s)}\right) = 0 \quad (30)$$

where $p(s)$ and $\tilde{p}(s)$ are as in (26). From the root locus technique, [18], the zeros of $p_2(s)$ converge to those of $p(s)$, namely, $-(\mu+\omega) < 0$ and $-\frac{1}{2}\big(2\mu+\sigma+\gamma \pm |\sigma-\gamma|\big) < 0$, as $\beta \to 0$. As a result, the eigenvalues of the linearized system (27) about $x_2^*$ are all stable from the continuity of the root locus of the characteristic equation for $\beta > 0$ not exceeding some sufficiently small threshold value for any given values of the remaining parameters of (8)-(10). Equivalently, that property holds if

$$\left\|\frac{\tilde{p}(s)}{p(s)}\right\|_\infty \triangleq \max_{\omega\in R_{0+}}\left\{\left|\frac{\tilde{p}(j\omega)}{p(j\omega)}\right|\right\} < \frac{1}{\beta} \quad (32)$$

This follows since $p(s)$ being a Hurwitz polynomial implies that $p_2(s)$ is Hurwitz if $\left|\beta \tilde{p}(j\omega)\right| < \left|p(j\omega)\right|$ $\forall \omega \in \mathbb{R}_{0+}$ from Rouché theorem of number of zeros within a closed set applied to the complex half-plane $\text{Re}\{s\} < 0$, [19].

Finally, note that the global Lyapunov stability is automatically guaranteed for the SEIR model (1)-(4) since the total population is assumed to be constant for all time. ***

**Remark 4**. The equilibrium points in the vaccination-free case are not suitable since one of them is concerned with the whole population being susceptible while the other one is concerned with not all the population being asymptotically converging to the removed-by-immunity, in general. Therefore, a suitable vaccination strategy is necessary to avoid the persistence of the disease in the population. Also, note that if $\beta < \dfrac{(\mu+\gamma)(\mu+\sigma)}{\sigma}$ then the equilibrium point $x_2^*$ does not exists, i.e., $x_1^*$ is the unique equilibrium point and it is stable. Otherwise, i.e., if $\beta \geq \dfrac{(\mu+\gamma)(\mu+\sigma)}{\sigma}$ then both equilibrium points exist and $x_1^*$ is an unstable attractor while the stability around $x_2^*$ depends on the fulfillment of the condition in Theorem 2 (ii). ***

## 3. VACCINATION STRATEGY

A control strategy may be defined in several ways involving the vaccination function which is really the manipulated variable. The control goal is decreasing appropriately the numbers of susceptible, infected and infectious while increasing the removed-by-immunity population. In this sense, two alternative vaccination strategies are proven to be appropriate from a health point of view. The following result is concerned with this matter.

**Theorem 3 (Control law 1)**. Assume that $\min\{S(0), E(0), I(0), R(0)\} \geq 0$, $\max\{S(0), E(0), I(0), R(0)\} \leq N$ and the feedback control with its associated vaccination strategy are as follows:

$$u(t) = -gS(t) \quad ; \quad g > \min\{\sigma, \gamma\} + \dfrac{k_0 \beta S(0)}{N} > 0 \tag{33}$$

$$V(t) = \dfrac{1}{\mu N}\left[\omega R(t) + \left(g - \dfrac{\beta I(t)}{N}\right)S(t) + \mu N\right] \tag{34}$$

for some $k_0 \in \mathbb{R}_+$. Then, the following results are derived:

**(i)** $S(t)$, $E(t)$, $I(t)$ and $R(t)$ are bounded $\forall t \in \mathbb{R}_{0+}$,

**(ii)** $S(t) \geq 0$, $E(t) \geq 0$ and $I(t) \geq 0$ $\forall t \in \mathbb{R}_{0+}$,

**(iii)** the whole population becomes asymptotically removed-by-immunity at an exponential rate,

and

**(iv)** there exist values for the control parameter $g$ guaranteeing that $V(t) \geq 0$ $\forall t \in \mathbb{R}_{0+}$.

*Proof*:

**(i)** Rewrite (1) in the equivalent form:

$$\dot{S}(t) = -\mu S(t) + u(t) \tag{35}$$

with an auxiliary control being defined as follows:

$$u(t) = \omega R(t) - \frac{\beta}{N} S(t)I(t) + \mu N (1 - V(t)) = -gS(t) \tag{36}$$

through the vaccination function given by (34). One gets from (35) and (36) that,

$$\dot{S}(t) = -(\mu + g)S(t) \Rightarrow S(t) = e^{-(\mu+g)t} S(0) \to S(\infty) = 0 \quad as \quad t \to \infty \tag{37}$$

since $g > 0$ and then $S(t)$ is bounded $\forall t \in \mathbb{R}_{0+}$. From (2)-(3) and (37), it follows that:

$$\dot{z}(t) = \underbrace{\begin{bmatrix} -(\mu+\sigma) & 0 \\ \sigma & -(\mu+\gamma) \end{bmatrix}}_{F_0} z(t) + \underbrace{\begin{bmatrix} 0 & \frac{\beta S(0)}{N} e^{-(\mu+g)t} \\ 0 & 0 \end{bmatrix}}_{\tilde{F}(t)} z(t) \tag{38}$$

where $z(t) = \begin{bmatrix} E(t) & I(t) \end{bmatrix}^T$. Note that the eigenvalues of $F_0$ are strictly negative, namely, $\lambda_1 = -(\mu+\sigma)$ and $\lambda_2 = -(\mu+\gamma)$. From (38), one gets that:

$$\|z(t)\| \leq k_0 e^{-\lambda_0 t} \|z(0)\| + \int_0^t k_0 e^{-\lambda_0(t-\tau)} \|\tilde{F}(\tau) z(\tau)\| d\tau \tag{39}$$

where $\lambda_0 = \min\{|\lambda_1|, |\lambda_2|\} = \mu + \min\{\sigma, \gamma\} > 0$. By direct calculations, it follows that:

$$\|z(t)\| \leq k_0 e^{-\lambda_0 t} \|z(0)\| + \frac{k_0 \beta S(0)}{N} e^{-\lambda_0 t} \int_0^t e^{(\lambda_0 - \mu - g)\tau} \|z(\tau)\| d\tau$$

$$\leq k_0 e^{-\lambda_0 t} \|z(0)\| + \frac{k_0 \beta S(0) e^{-\lambda_0 t}}{N(\lambda_0 - \mu - g)} \left( e^{(\lambda_0 - \mu - g)t} - 1 \right) \max_{0 \leq \tau \leq t} \{\|z(\tau)\|\} \leq k_0 \|z(0)\| + \frac{k_0 \beta S(0)}{N(g - \min\{\sigma, \gamma\})} \max_{0 \leq \tau \leq t} \{\|z(\tau)\|\}$$

$$\tag{40}$$

where the condition for the control parameter $g$ in (33) has been taking into account. From the fact that the right hand side of (40) is a monotonically increasing function, it follows that:

$$\max_{0 \leq \tau \leq t} \{\|z(\tau)\|\} \leq k_0 \frac{N(g - \min\{\sigma, \gamma\})}{N(g - \min\{\sigma, \gamma\}) - k_0 \beta S(0)} \|z(0)\| < \infty \tag{41}$$

provided that $g$ satisfies the condition in (33). Then, $E(t)$ and $I(t)$ are bounded $\forall t \in \mathbb{R}_{0+}$. Finally, $R(t)$ is bounded $\forall t \in \mathbb{R}_{0+}$ from Assertion 1.

**(ii)** From (37), $S(t) \geq 0$ $\forall t \in \mathbb{R}_{0+}$. Moreover, note that $z: \mathbb{R}_{0+} \to \mathbb{R}^2$ is continuous since it satisfies

(38) and consider $z(0) \geq 0$, i.e., both $E(0) \geq 0$ and $I(0) \geq 0$. Then, the results $E(t) \geq 0$ and $I(t) \geq 0$ $\forall t \in \mathbb{R}_{0+}$ are proved by contradiction. For such a purpose, consider the first time instant $t \in \mathbb{R}_{0+}$ such that $z(t) < 0$, i.e., $E(t) < 0$ and/or $I(t) < 0$. Then, $z(\tau) \geq 0$ $\forall \tau \in [0, t)$ and

$$z(t) = e^{F_0 t} z(0) + \int_0^t e^{F_0(t-\tau)} \tilde{F}(\tau) z(\tau) d\tau \geq 0 \quad (42)$$

since $e^{F_0 t} > 0$ $\forall t \in \mathbb{R}_{0+}$, due to $F_0$ is a Metzler matrix, and $e^{F_0(t-\tau)} \tilde{F}(\tau) z(\tau) \geq 0$, due to $\tilde{F}(\tau) \in \mathbb{R}_{0+}^{2 \times 2}$ by construction. Thus, (42) contradicts $z(t) < 0$ and there is not $t \in \mathbb{R}_{0+}$ such that $z(t) < 0$. As a consequence, $E(t) \geq 0$ and $I(t) \geq 0$ $\forall t \in \mathbb{R}_{0+}$.

**(iii)** From (2) and (37), it follows that:

$$\dot{E}(t) = -(\mu + \sigma) E(t) + \frac{\beta S(0)}{N} e^{-(\mu+g)t} I(t) \quad (43)$$

what leads to:

$$E(t) = e^{-(\mu+\sigma)t} E(0) + \frac{\beta S(0)}{N} e^{-(\mu+\sigma)t} \int_0^t e^{-(g-\sigma)\tau} I(\tau) d\tau$$

$$\leq \begin{cases} e^{-(\mu+\sigma)t} E(0) + \frac{\beta S(0) M_I}{N} \frac{e^{-(\mu+g)t} - e^{-(\mu+\sigma)t}}{\sigma - g} & \text{if } g \neq \sigma \\ e^{-(\mu+\sigma)t} E(0) + \frac{\beta S(0) M_I}{N} t e^{-(\mu+\sigma)t} & \text{if } g = \sigma \end{cases} \quad (44)$$

where $M_I = \max_{0 \leq \tau \leq t} \{I(\tau)\} < \infty$. Such an upper-bound for $E(t)$ tends exponentially fast to zero as $t \to \infty$ at a rate of at most $\mu + \min\{\sigma, g\}$ for any $g > 0$. From (3), it follows:

$$I(t) = e^{-(\mu+\gamma)t} I(0) + \sigma \int_0^t e^{-(\mu+\gamma)(t-\tau)} E(\tau) d\tau \quad (45)$$

By introducing (44) in (45), one gets that:

$$I(t) \leq e^{-(\mu+\gamma)t} I(0) + \frac{\sigma E(0)}{\gamma - \sigma} \left( e^{-(\mu+\sigma)t} - e^{-(\mu+\gamma)t} \right) + \frac{\sigma \beta S(0) M_I}{N(\sigma - g)} \left( \frac{e^{-(\mu+g)t} - e^{-(\mu+\gamma)t}}{\gamma - g} - \frac{e^{-(\mu+\sigma)t} - e^{-(\mu+\gamma)t}}{\gamma - \sigma} \right) \quad (46)$$

in case that $g \neq \sigma$, $g \neq \gamma$ and $\sigma \neq \gamma$. For another relations between $\sigma$, $\gamma$ and $g$, analogy calculations lead to alternative expressions to (46). In any case, such expressions show that $I(t)$ tends exponentially fast to zero as $t \to \infty$ for any $g > 0$. Finally, $R(t) = (N - S(t) - E(t) - I(t))$ tends exponentially fast to $N$ as $t \to \infty$ for any $g > 0$.

**(iv)** From (37), $S(t) \in [0, N]$ $\forall t \in \mathbb{R}_{0+}$. Moreover, both $I(t)$ and $R(t)$ are bounded $\forall t \in \mathbb{R}_{0+}$ as it has been proven in part (i) of this theorem. Then, the existence of values for the control parameter $g$ guaranteeing that $V(t) \geq 0$ $\forall t \in \mathbb{R}_{0+}$ is derived from (34). ∗∗∗

**Remark 5.**
**(i)** The following notation has been used in the part (ii) of Theorem 2:

- $M \in \mathbb{R}_{0+}^{n \times n}$ and $M > 0$ denote a positive real matrix in the usual sense that all its entries are nonnegative,
- $v < 0$ denotes a vector with at least one component being negative.

**(ii)** Note that if the control parameter $g$ is such that $g > \min\{\sigma, \gamma\} + k_0 \beta$ the same results of Theorem 3 are satisfied irrespective of the initial condition for the susceptible population.  ***

An alternative vaccination strategy is dealt with in below:

**Theorem 4 (Control law 2)**. Assume that the vaccination strategy is as follows:

$$V(t) = \begin{cases} \dfrac{1}{\mu N}\big(g_1 N - g\, R(t) - \gamma\, I(t)\big) & for\ 0 \le t < t_s \\ \dfrac{1}{\mu N}\big(\mu N + \omega\, R(t)\big) & for\ t \ge t_s \end{cases} \quad (47)$$

with $g > -(\mu + \omega)$ and $g_1 \ge \mu + \omega + g > 0$, where $t_s$ is the time instant, if it exists, at which the susceptible population becomes zero, i.e., $S(t_s) = 0$. The above vaccination strategy implies the following results:

**(i)** $\min\{S(t), E(t), I(t), R(t)\} \ge 0$ and $\max\{S(t), E(t), I(t), R(t)\} \le N$ $\forall t \in \mathbb{R}_{0+}$ irrespective of the existence or not of the time instant $t_s$ provided that $\min\{S(0), E(0), I(0), R(0)\} \ge 0$.

**(ii)** If the time instant $t_s$ does not exist, i.e., if the susceptible population does not become zero at any time, then:

$$R(\infty) = \frac{g_1 N}{\mu + \omega + g} \ ; \ S(\infty) + E(\infty) + I(\infty) = \frac{(\mu + \omega + g - g_1)N}{\mu + \omega + g} \quad (48)$$

irrespective of the initial conditions. In particular, if $g_1 = \mu + \omega + g$ then $R(\infty) = N$ and $S(\infty) + E(\infty) + I(\infty) = 0$, i.e., the removed-by-immunity population equalizes asymptotically the total population at exponential rate while the sum of the infected, infectious and susceptible populations converges asymptotically to zero at exponential decay rate.

**(iii)** If $t_s$ exists then $R(\infty) = N$ and $S(\infty) = E(\infty) = I(\infty) = 0$.

**(iv)** In the case that $g_1 = \mu + \omega + g$ and provided that $\min\{S(0), E(0), I(0), R(0)\} \ge 0$, the vaccination function is nonnegative for all time if $g_1 \ge \gamma$ or if $g \ge 0$ and $\gamma = \mu + \omega$.

*Proof*:
**(i)** Suppose that $t_s$ does not exist (or $t_s \to \infty$), then $S(t) > 0$ $\forall t \in \mathbb{R}_{0+}$ provided that $S(0) > 0$ from the definition of $t_s$. By introducing the vaccination control law (47) in (4) one obtains:

$$\dot{R}(t) = -(\mu + \omega + g)R(t) + g_1 N \quad (49)$$

Then, $R(t) = e^{-(\mu+\omega+g)t}\left(R(0) + g_1 N \int_0^t e^{(\mu+\omega+g)\tau} d\tau\right) > 0$ $\forall t \in \mathbb{R}_{0+}$ if $R(0) > 0$. Finally, $E(t) > 0$ and $I(t) > 0$ $\forall t \in \mathbb{R}_{0+}$ provided that $E(0) > 0$ and $I(0) > 0$ by using similar arguments that those used in Theorem 3 (ii). Otherwise, i.e., if $t_s$ exists then $S(t) > 0$, $E(t) > 0$, $I(t) > 0$ and $R(t) > 0$ $\forall t \in [0, t_s)$ provided that $\min\{S(0), E(0), I(0), R(0)\} \geq 0$ from similar arguments to those used in the case that $t_s$ does not exist (or $t_s \to \infty$). Moreover, $S(t) = 0$ $\forall t \geq t_s$ from the definition of the control law (47) and then $E(t) = E(t_s) e^{-(\mu+\sigma)(t-t_s)} \geq 0$ $\forall t \geq t_s$ from (2) by taking into account that $E(t) > 0$ $\forall t \in [0, t_s)$ and it is a continuous function. Furthermore, it follows that:

$$I(t) = e^{-(\mu+\gamma)(t-t_s)} I(t_s) + \frac{\sigma E(t_s)}{\gamma - \sigma}\left(e^{-(\mu+\sigma)(t-t_s)} - e^{-(\mu+\gamma)(t-t_s)}\right) \quad \forall t \geq t_s \quad \text{if} \quad \gamma \neq \sigma$$

or (50)

$$I(t) = e^{-(\mu+\gamma)(t-t_s)} I(t_s) + \sigma E(t_s) e^{-(\mu+\gamma)(t-t_s)}(t-t_s) \quad \forall t \geq t_s \quad \text{if} \quad \gamma = \sigma$$

from (3). Then, $I(t) \geq 0$ $\forall t \geq t_s$ since $I(t) > 0$ $\forall t \in [0, t_s)$ and it is a continuous function. By introducing the vaccination function (47) in (4), one obtains that $\dot{R}(t) = -\mu R(t) + \gamma I(t) + \mu N$ $\forall t \geq t_s$, so that:

$$R(t) = e^{-\mu(t-t_s)} R(t_s) + \int_0^t e^{-\mu(t-\tau)}(\gamma I(\tau) + \mu N) d\tau \geq 0 \quad \forall t \geq t_s \quad (51)$$

since $I(t) \geq 0$ $\forall t \geq t_s$ and $R(t) > 0$ $\forall t \in [0, t_s)$ and it is a continuous function. In summary, $\min\{S(t), E(t), I(t), R(t)\} \geq 0$ $\forall t \in \mathbb{R}_{0+}$ irrespective of the existence of the time instant $t_s$ provided that $\min\{S(0), E(0), I(0), R(0)\} \geq 0$. As a result, $\max\{S(t), E(t), I(t), R(t)\} \leq N$ $\forall t \in \mathbb{R}_{0+}$ from Assertion 1.

(ii) Suppose that $t_s$ does not exist (or $t_s \to \infty$), then by introducing the vaccination control law (47) in (4) one obtains (49) and, as a consequence, $R(t) \to R(\infty) \triangleq \lim_{t\to\infty}\{R(t)\} = \frac{g_1 N}{\mu + \omega + g}$ as $t \to \infty$ since $g > -(\mu+\omega)$. Thus, from Assertion 1 it follows that $S(\infty) + E(\infty) + I(\infty) = \frac{(\mu+\omega+g-g_1)N}{\mu+\omega+g}$. In the particular case that $g_1 = \mu+\omega+g$, $R(\infty) = N$ and $S(\infty) + E(\infty) + I(\infty) = 0$ from the above expressions.

(iii) If $t_s$ exists, from (1) and (47), it is deduced that:

$$\dot{S}(t) = \omega R(t) + \mu N(1 - V(t)) = 0 \quad \forall t \geq t_s \quad (52)$$

since both $S(t) = 0$ and $\dot{S}(t) = 0$ $\forall t \geq t_s$. Furthermore, $E(t) = E(t_s) e^{-(\mu+\sigma)(t-t_s)}$ $\forall t \geq t_s$ from (2) and, as a consequence, $E(t) \to E(\infty) = 0$ as $t \to \infty$ since $\mu + \sigma > 0$. Moreover, the expression (50) for $I(t)$ is deduced from (3) so that $I(t) \to I(\infty) = 0$ as $t \to \infty$ since $\mu + \gamma > 0$. Finally, from Assertion 1 $R(t) \to R(\infty) = N$ as $t \to \infty$ is obtained.

**(iv)** In the case that $g_1 = \mu + \omega + g$, the vaccination function is nonnegative from (47) for all time if $g_1 \geq \gamma$ since:

$$g_1 N - gR - \gamma I = g_1(R + S + E + I) - gR - \gamma I = (g_1 - g)R + (g_1 - \gamma)I + g_1(S + E)$$
$$= (\mu + \omega)R + (g_1 - \gamma)I + g_1(S + E) \geq 0 \qquad (53)$$

where the fact that $S(t) \geq 0$, $E(t) \geq 0$, $I(t) \geq 0$ and $R(t) \geq 0$ $\forall t \in \mathbb{R}_{0+}$ has been used. Furthermore, if the constraint $g_1 \geq \gamma$ is changed to $g \geq 0$ and $\gamma = \mu + \omega$ then $g_1 - g = \gamma = \mu + \omega$, equivalently $g_1 - \gamma = g \geq 0$, and the vaccination function is also nonnegative for all time since the above expression becomes:

$$g_1 N - gR - \gamma I = (\mu + \omega)R + gI + g_1(S + E) \geq 0 \qquad (54)$$
***

**Remark 6.** In the particular case that $g_1 = \mu + \omega + g$, one gets from (47) for all time that:

$$I(t) \leq \frac{(g_1 - g)N}{\gamma} = \frac{(\mu + \omega)N}{\gamma} \leq \frac{g_1 N - gR(t)}{\gamma} \Rightarrow V(t) \geq 0 \qquad (55)$$

since $R(t) \leq N$ for all time. Such a result is guaranteed for arbitrary initial conditions of (1) from Assertion 1 if $g_1 - g = \gamma = \mu + \omega$. ***

**Corollary 1.** The vaccination strategy of Theorem 4 is nonnegative for all time if

$$\min\left\{1, \frac{\sigma\beta}{(\mu + \sigma)(\mu + \gamma)}\right\} \leq \frac{\mu + \omega}{\gamma}.$$

**Proof**: From (3), $I(t) = \frac{\sigma E(t)}{\mu + \gamma}$ with $E(t) \neq 0$ is a potential maximum of the infectious population since it makes $\dot{I}(t) = 0$. Thus, $I(t) \leq I_{max} \triangleq \max\{I(t) : t \geq 0\} \leq \frac{\sigma E_{max}}{\mu + \gamma}$, where $E_{max} \triangleq \max\{E(t) : t \geq 0\}$ is reached for $\dot{E}(t) = 0$ in (2) so that $(\mu + \sigma)E_{max} \leq \beta \frac{S_{max} I_{max}}{N}$, where $S_{max} \triangleq \max\{S(t) : t \geq 0\}$. Combining the two relations and using Assertion 1 one has for all time:

$$I(t) \leq I_{max} \leq \frac{\sigma E_{max}}{\mu + \gamma} \leq \frac{\sigma\beta S_{max} I_{max}}{(\mu + \sigma)(\mu + \gamma)N} \leq \frac{\sigma\beta N}{(\mu + \sigma)(\mu + \gamma)} \Rightarrow I(t) \leq I_{max} \leq \min\left\{1, \frac{\sigma\beta}{(\mu + \sigma)(\mu + \gamma)}\right\} N$$
(56)

Then, the proof is derived by combining (55) and (56). ***

## 4. SIMULATION RESULTS

A pair of examples for illustrating the theoretical results presented in the paper are described in this section.

### 4.1. Example 1: A measles infection

An example based on an outbreak of measles in a total population of $N = 10^6$ habitants is considered. Such an epidemic can be described by the SEIR model (1)-(4) with the parameter values: $\mu = 5.48 \times 10^{-5}$ per day ($p.\,d.$), $\beta = 3.288\ p.\,d.$, $\omega = 0\ p.d.$ $\sigma = 9.82 \times 10^{-2}\ p.d.$ and $\gamma = 0.274\ p.d.$.

Such values are commonly used in the literature, [14]. The initial conditions for the individual populations are given by: $S(0) = 9.8 \times 10^5$, $E(0) = 1.5 \times 10^4$, $I(0) = 5000$ and $R(0) = 0$.

The time evolution of the system in the free-vaccination case converges to the equilibrium point $x_2^* = \begin{bmatrix} S_2^* & I_2^* & R_2^* \end{bmatrix}^T$ with the individual populations obtained from (23)-(25). The percentages of susceptible, infected, infectious and removed-by-immunity populations with respect to the total population at such an equilibrium point are $S_2^*(\%) = 8.34$, $E_2^*(\%) = 100 - S_2^*(\%) - I_2^*(\%) - R_2^*(\%) = 0.051$, $I_2^*(\%) = 0.019$ and $R_2^*(\%) = 91.59$, respectively. This is an inappropriate equilibrium point since there would be an appreciable number of susceptible individuals. i.e, the infectious disease would not be eradicated from the population. As a consequence, a control action could be applied in order to eradicate the disease.

*4.1.1. Epidemic evolution with the vaccination control law 1*

The vaccination strategy associated with the control law 1 defined by (33)-(34) is considered. The control parameter $g = 0.25$ is applied. The time evolution of the respective populations (percentage respect the total population) is displayed in Figure 1. It can be seen that the model tends to a suitable equilibrium point as time goes to infinity since all the population becomes asymptotically removed-by-immunity. i.e. the infection would be eradicated from the population with such a vaccination practice in a relative short time period, approximately 50 days. Figure 2 displays the time evolution of the vaccination effort $\mu N V(t)$ (in terms of percentage with respect to the total population) to be applied to eradicate the disease.

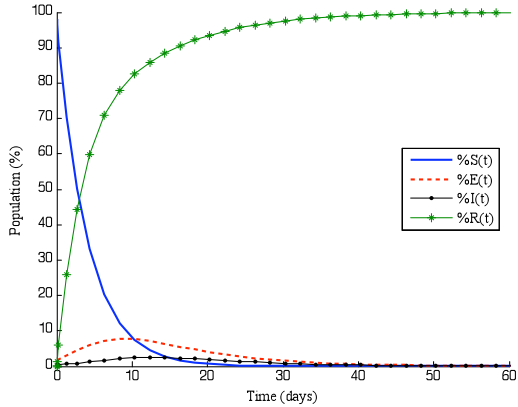 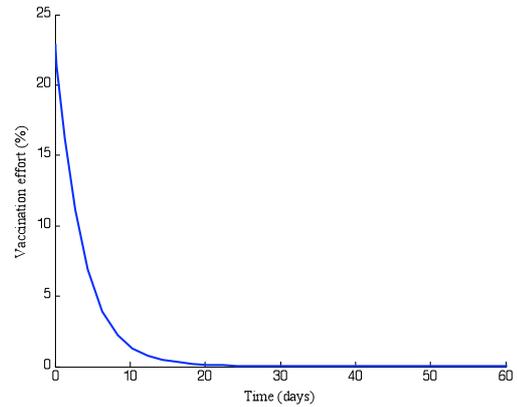

Fig. 1. Time evolution of the partial populations with the vaccination control law 1.

Fig. 2. Time evolution of the vaccination effort $\mu N V(t)$ with the control law 1.

*4.1.2. Epidemic evolution with the vaccination control law 2*

A vaccination strategy based on the control law 2 defined by (47) with the parameters $g = 0.0999$ and $g_1 = \mu + \omega + g = 0.1$ is applied. The time evolution of the partial populations (percentage respect the total population) is displayed in Figure 3. Again a suitable equilibrium point is reached as time tends to infinity since all the population becomes asymptotically removed-by-immunity. As a result, the infection is eradicated from the population in 50 days, approximately. Figure 4 displays the time evolution of the vaccination effort $\mu NV(t)$ (in terms of percentage with respect to the total population) to be applied to eradicate the disease.

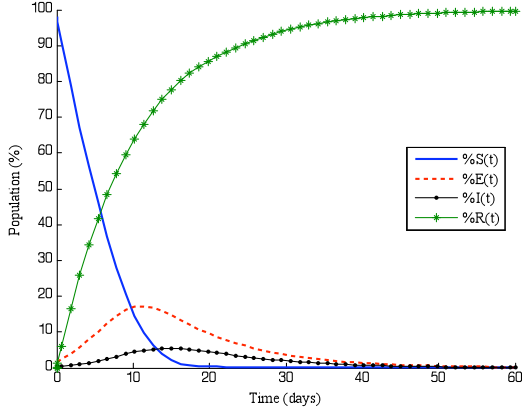 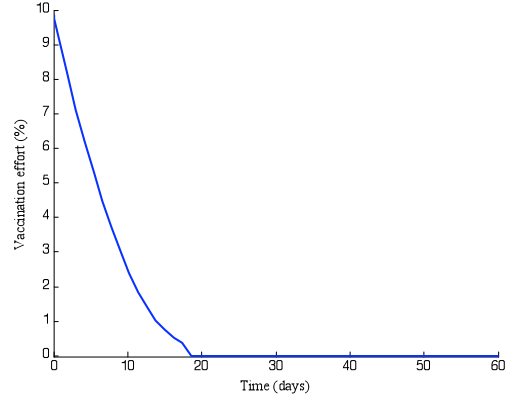

Fig. 3. Time evolution of the partial populations with the vaccination control law 2.

Fig. 4. Time evolution of the vaccination effort $\mu NV(t)$ with the control law 2.

*4.1.3. Design of a vaccination campaign*

The time evolution of the vaccination effort may be used to prescribe a suitable vaccination campaign in the population. For such a purpose, note that all the population becomes removed-by-immunity at the equilibrium point with both vaccination laws (see figures 1 and 3) while 8.4% of the population did not reach such a status at the equilibrium point of the vaccination-free case, approximately. The later percentage determines the population ($n = 84000$ habitants) to be vaccinated in order to eradicate the disease. In view of figures 2 and 4, the vaccination effort depends on time through a uniformly decreasing function. Then, the vaccination cadence will not be uniform during the disease evolution (approximately 50 days) but it presents a maximum along the first day and will be uniformly decreasing to a zero value at the 20$^{th}$ day, approximately. In this sense, the vaccination cadence associated to any of the control laws may be given by:

$$v(t) = \frac{\mu NV(t)}{f} \quad \text{vaccines per day,} \tag{57}$$

where the normalization factor $f$ is added *to relate as recommendation the vaccination effort with the number of vaccines to be applied*. For the purpose of calculating $f$, the contribution of the vaccination effort to the removed-by-immunity population is given by $\mu N \int_0^{50} V(t)\, dt$ from (4), so that:

$$\frac{\mu N}{f} \int_0^{50} V(t)\, dt = n \tag{58}$$

provided that each individual is vaccinated at most once. From figures 2 and 4, one can evaluate $\mu N \int_0^{50} V(t)\, dt = 8.385 \times 10^5$ for the control law 1 (then $f = 9.98$) and $\mu N \int_0^{50} V(t)\, dt = 6.65 \times 10^5$ for the control law 2 (then $f = 7.92$). The figures 5 and 6 display the vaccination cadence per day, i.e. the number of vaccines to be applied each day, according to (57) for both vaccination strategies.

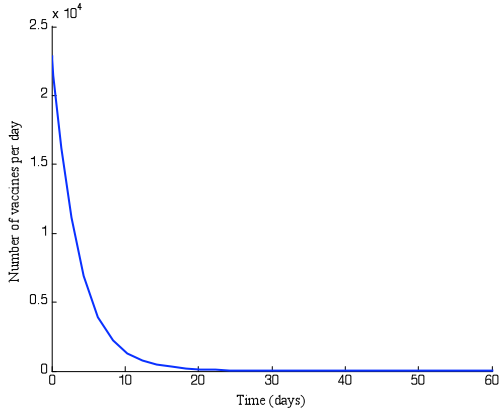 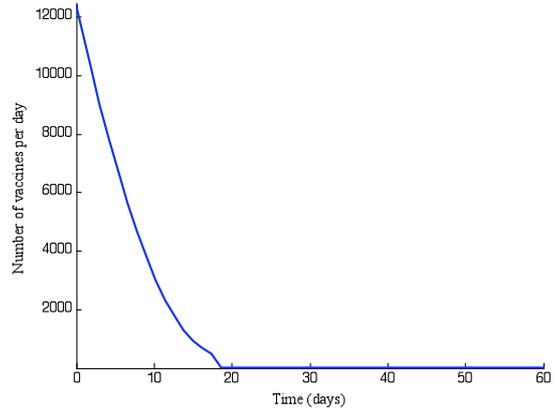

Fig. 5. Vaccination cadence per day for control law 1.   Fig. 6. Vaccination cadence per day for control law 2.

*4.1.4. Control parameters influence in the infection time evolution*

In one hand, the disease propagation experiments a peak in the infectious population (2.5% of the total population) at the 12$^{th}$ day with the application of the vaccination control law 1 as it can be seen from Figure 1. On the other hand, the use of the vaccination control law 2 makes the disease reaches a peak at the 15$^{th}$ day when the 5.4% of the total population is infectious as it can be seen from Figure 3. Both specifications (the maximum of the infection and the peak time) depend on the values of the control parameters ($g$ for the control law 1 and, $g$ and $g_1$ for the control law 2). In this sense, Figures 7(a) and 7(b) display the evolution of the infectious population for different values of such control parameters for both vaccination strategies.

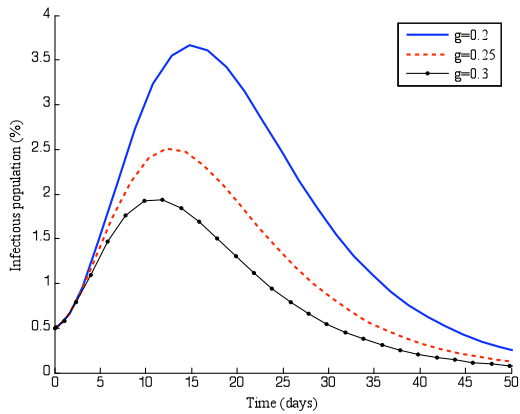 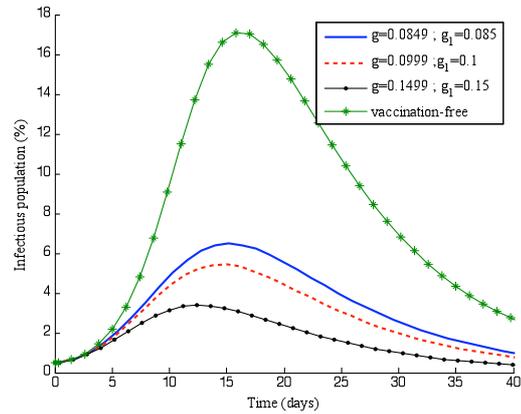

(a)                                                                   (b)

Figs. 7. Time evolution of the infectious population with different values of the control parameters for (a) the vaccination control law 1 and (b) the vaccination control law 2

*4.2. Example 2: An influenza infection*

An example based on an outbreak of influenza in a British boarding school in early 1978, [8], is used. Such an epidemic can be described by the SEIR mathematical model (1)-(4) with the parameter values $\frac{1}{\mu} = 70\ years = 25550\ days$, $\beta = 1.66\ per\ day$, $\frac{1}{\sigma} = \frac{1}{\gamma} = 2.2\ days$ and $\frac{1}{\omega} = 7\ days$ (case 1) or $\frac{1}{\omega} = 15\ days$ (case 2). A total population of $N = 1000\ boys$ is considered with the initial conditions for the individual populations given by $S(0) = 980\ boys$, $E(0) = 15\ boys$, $I(0) = 5\ boys$ and $R(0) = 0\ boys$.

Two sets of simulation results are presented to compare the evolution of the SEIR mathematical model populations in two different situations, namely: when no control actions are applied and if a control action is applied via a vaccination practice.

*4.2.1. Epidemic time evolution without vaccination*

The time evolution of the respective populations is displayed in Figure 8(a) (for the case that $\frac{1}{\omega} = 7\ days$) and Figure 8(b) (if $\frac{1}{\omega} = 15\ days$). The model tends to an inappropriate equilibrium point as time goes to infinity in both cases. Note that there would be susceptible, infected and infectious population for all time. i.e, the infectious disease would not be eradicated. As a consequence, a control action will have to be applied in order to eradicate the disease. Also, the number of removed-by-immunity individuals in the equilibrium point decreases if the loss of immunity rate $\omega$ increases, and conversely, as it can be deduced from such figures. At the same time, the number of infected and infectious individuals in the equilibrium point increases if the loss of immunity rate $\omega$ increases, and conversely. Finally, the susceptible population in the equilibrium point is practically the same for both values of $\omega$.

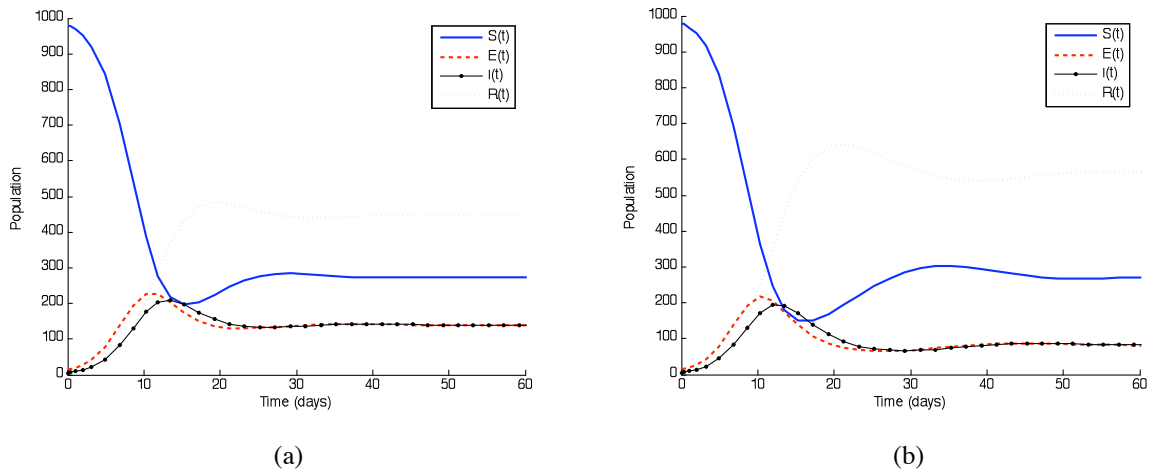

Figs. 8. Time evolution of the individual populations without vaccination for
(a) $\frac{1}{\omega} = 7\ days$ and (b) $\frac{1}{\omega} = 15\ days$

*4.2.1. Epidemic time evolution with the vaccination control law 1*

A feedback control law defined by (33)-(34) with the control parameter $g = 0.1$ is applied. The time evolution of the respective populations is displayed in Figure 9 so in the case that $1/\omega = 7\ days$ as $1/\omega = 15\ days$. Note that the use of such a control law makes the dynamics of the susceptible population be given by $\dot{S}(t) = -(\mu + g)S(t)$, i.e. its time evolution does not depend on the loss of immunity rate $\omega$. Then, the time evolution of the infected and infectious population is also independent of the loss of immunity rate from (2) and (3). Finally, the removed-by-immunity population does not depend on such a parameter from Assertion 1. In summary, the time evolution of all partial populations is independent of $\omega$ due to the vaccination control rule (33)-(34). Furthermore, it can be seen that the model tends to a suitable equilibrium point as time goes to infinity since all population become removed-by-immunity. i.e. the infection would be eradicated from the population with such a vaccination practice in a relative short time period, approximately 49 days (7 weeks). Figure 10 displays the time evolution of the vaccination effort $\mu N V(t)$ to be applied to eradicate the disease, for both values of $\omega$.

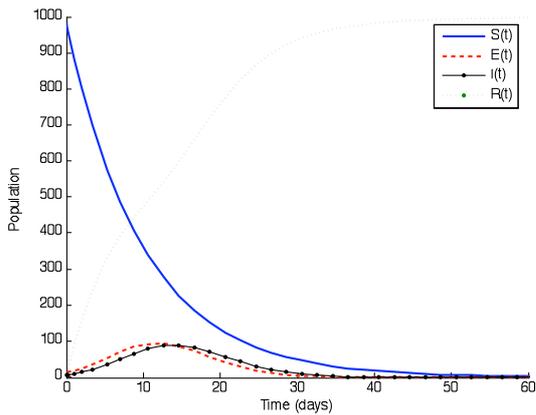 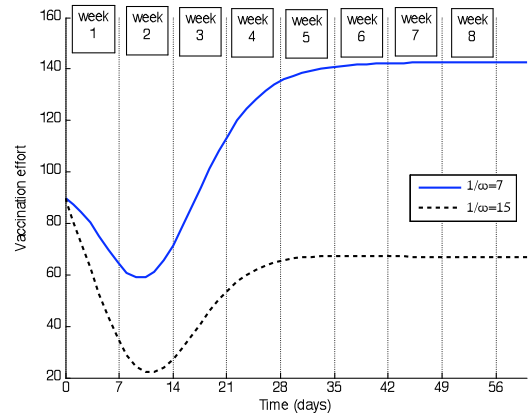

Fig. 9. Time evolution of the individual populations if the vaccination control law 1 with $g = 0.1$ is applied.

Fig. 10. Time evolution of the vaccination effort $\mu N V(t)$ for $1/\omega = 7\ days$ and for $1/\omega = 15\ days$.

The time evolution of the vaccination effort may be used to prescribe a suitable vaccination campaign in the population. For such a purpose, note that the removed-by-immunity population is $r_v = 993$ boys at the steady state day 49 with the application of such a vaccination law (see Figure 9) and it is $r_1 = 446$ boys, if $1/\omega = 7\ days$, or $r_2 = 562$ boys, if $1/\omega = 15\ days$, without a vaccination effort as it can be seen from figures 8(a) and 8(b), respectively. Then, the number of boys to be vaccinated is given by $n_1 = r_v - r_1 = 547$ if $1/\omega = 7\ days$ or by $n_2 = r_v - r_2 = 431$ if $1/\omega = 15\ days$, provided that each boy receives at most one vaccine. In view of Figure 10, the vaccination cadence is not uniform during the disease evolution (7 weeks) but it presents a minimum during the second week and it reaches a maximum during the last four weeks. In this sense, the vaccination cadence may be given by:

$$v_i(t) = \frac{\mu N V_i(t)}{f_i} \text{ vaccines per day} \tag{59}$$

for $i = 1$ if $1/\omega = 7\ days$ or $i = 2$ if $1/\omega = 15\ days$. Again, the normalization factor $f_i$ is added *to relate as recommendation the vaccination effort with the number of vaccines to be applied.* In this case, from Figure 10 $f_1 = 10$ and $f_2 = 6.3$ are obtained so that:

$$\frac{\mu N}{f_i} \int_0^{49} V_i(t)\, dt = n_i \tag{60}$$

for $i = 1, 2$. Figure 11 displays the vaccination cadence per day, i.e. the number of vaccines to be applied each day, for both values of $\omega$. However, *a vaccination cadence per week may be recommended to be used in practice*. Then, the vaccination cadence per day may be given by:

$$v_i(t) = \frac{\mu N}{T f_i} \int_{(j-1)T}^{jT} V_i(t)\, dt \quad for\ (j-1)T \leq t < jT \tag{61}$$

where $T = 7\ days$ (1 *week*) and $j = 1, 2, ..., 7$. In this way, the vaccines to be applied each week are uniformly distributed in the days of such a week. Figure 12 displays the vaccination cadence if this method for the vaccines distribution is used.

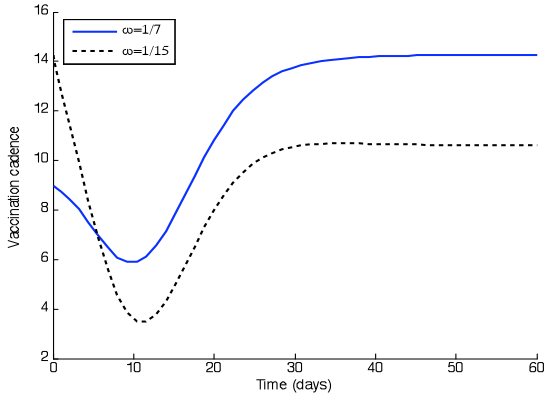 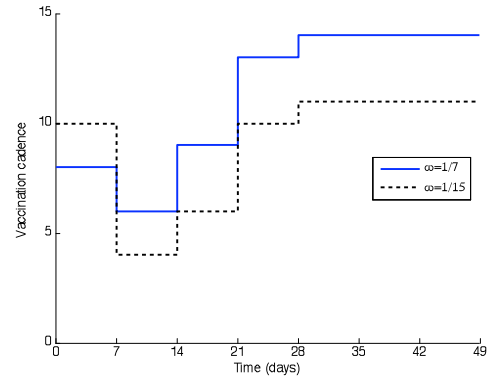

Fig. 11. Vaccination cadence per day for $1/\omega = 7\ days$ and for $1/\omega = 15\ days$.

Fig. 12. Weekly uniform vaccination cadence per day for $1/\omega = 7\ days$ and for $1/\omega = 15\ days$.

*4.2.2. Epidemic time evolution with the vaccination control law 2*

A feedback control law defined by (47) with the control parameters $g = -0.015$ and $g_1 = \mu + \omega + g$ (i.e., $g_1 = 0.1297$ for $1/\omega = 7\ days$ and $g_1 = 0.0517$ for $1/\omega = 15\ days$) is applied. The time evolution of the respective populations is displayed in figures 13.

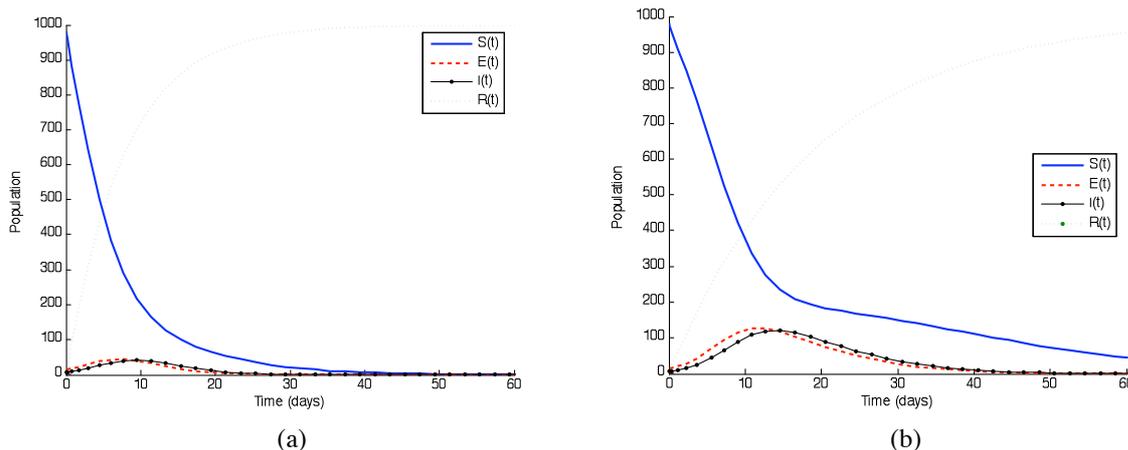

(a)                                             (b)

Figs. 13. Time evolution of the individual populations with the vaccination control law 2 with the parameters (a) $g = -0.015$ and $g_1 = 0.1297$ if $1/\omega = 7\ days$ and (b) $g = -0.015$ and $g_1 = 0.0517$ if $1/\omega = 15\ days$

Finally, Figure 14 displays the time evolution of the vaccination effort $\mu N V(t)$ to be applied to eradicate the disease, for both values of $\omega$. By using this result, a similar study to that carried out with the control law 1 could be done to establish an appropriated vaccination campaign through the population.

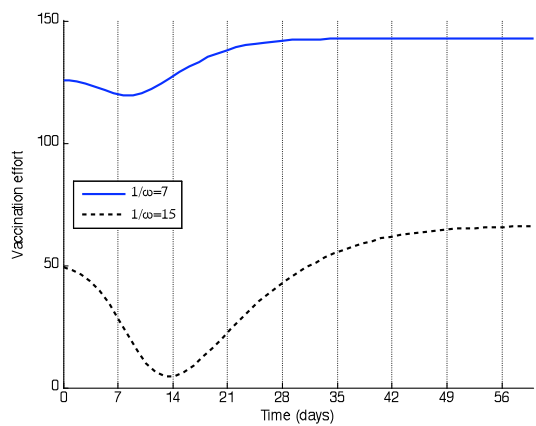

Fig. 14. Time evolution of the vaccination function for $1/\omega = 7\ days$ and for $1/\omega = 15\ days$.

## ACKNOWLEDGEMENTS

The authors are very grateful to the Spanish Ministry of Education by its partial support of this work through Project DPI2009-07197 and to the Basque Government by its support through Grants IT378-10 and SAIOTEK SPE07UN04.